\documentclass[twoside,reqno,A4]{amsart}

\numberwithin{equation}{section}

\newcommand{\qdn}{\hspace*{-1.5mm}}
\newcommand{\qqdn}{\hspace*{-2.5mm}}

\newcommand{\+}{&\qqdn}%

\newcommand{\mb}[1]{\mathbb{#1}}

\newcommand{\ffnk}[4]{\left[\qdn\ba{#1}#3\\[2mm]#4\ea{\!;\:#2}\right]}

\newcommand{\nnm}{\nonumber}
\newcommand{\be}{\begin{equation}}
\newcommand{\ee}{\end{equation}}
\newcommand{\ba}{\begin{array}}
\newcommand{\ea}{\end{array}}
\newcommand{\bmn}{\begin{eqnarray}}
\newcommand{\emn}{\end{eqnarray}}
\newcommand{\bnm}{\begin{eqnarray*}}
\newcommand{\enm}{\end{eqnarray*}}
\newcommand{\bln}{\begin{subequations}}
\newcommand{\eln}{\end{subequations}}

\newcommand{\bbtm}[4]{\bibitem{kn:#1}{#2,}~\emph{#3,}~{#4.}}
\newcommand{\cito}[1]{\cite{kn:#1}}

\makeatletter
\def\section{\@startsection {section}{1}{\z@}%
{-3.5ex \@plus -1ex \@minus -.2ex}%
{2.3ex \@plus.2ex}%
{\sectionformat}}
\def\sectionformat{\normalfont\large\sc}
\makeatother
\begin{document} 
\title{New Laplace transforms for the generalized \\ hypergeometric functions $_2F_2$ and $_3F_3$}
\dedicatory{\textsc{\large Xiaoxia Wang $^{a, *}$, Arjun K. Rathie $^b$}\\[1mm]
$^a$ Department of Mathematics, Shanghai University, Shanghai, 200444, P. R. China;\\[1mm]
$^b$ Department of Mathematics, School of Mathematical and Physical Sciences, Central University of Kerala, \\[1mm]
Riverside Transit Campus, Padennakkad P. O. Nileshwar Kasaragod-671 328, Kerala State, India.}

\thanks{$^{*}$ Corresponding author. \\
E-mail addresses: xiaoxiawang@shu.edu.cn (X. Wang), akrathie@rediffmail.com (A. K. Rathie) }


\maketitle\thispagestyle{empty}
\markboth{X. Wang and A.K. Rathie}
{New Laplace transforms for the generalized hypergeometric functions $_2F_2$ and $_3F_3$}

\begin{center}\parbox{120mm}{

Motivated by the new Laplace transforms for the Kummer's confluent hypergeometric functions $_1F_1$ obtained recently by Kim et al.
[Math $\&$ Comput. Modelling, 55 (2012), pp. 1068--1071], the authors aim is to establish so far unknown Laplace transforms of rather
general case of generalized hypergeometric functions $_2F_2(x)$ and $_3F_3(x)$ by employing extensions of classical summation theorems for the series
$_2F_1$ and $_3F_2$ obtained recently by Kim et al. [Int. J. Math. Math. Sci., 309503, 26 pages, 2010].
Certain known results obtained earlier by Kim et al. follow cases of our main findings. \\[4mm]

\textbf{Keywords:} Laplace transform;
                  Watson's summation theorem;
          Dixon's summation theorem;
          Whipple's summation theorem;
          extension summation theorem.\\[4mm]
\emph{2000 Mathematics Subject Classification}:  {Primary 33C20, Secondary 33C05, 33B20}.}
\end{center}

\section{Introduction and results required}
We begin with the definition of generalized hypergeometric function with $p$ numerator parameters and $q$ denominator parameters
($p$ and $q$ being nonnegative integers) by means of the following series \cite{kn:rainville, kn:slater}:
\bmn\label{def}
{_p F_q}
\ffnk{ccc}{z}{a_1,\+\cdots,\+a_p}
              {b_1,\+\cdots,\+b_q}
=\sum_{n=0}^{\infty}
\frac{(a_1)_n\cdots(a_p)_n}
     {(b_1)_n\cdots(b_q)_n}
 \cdot\frac{z^n}{n!}
\emn
whenever this series converges and elsewhere by analytic continuation.
Also $(\cdot)_n$ denotes the well known Pochhammer symbol (or the shifted factorial) defined by
\bmn
(a)_n
=\begin{cases}
1,\+n=0\\[2mm]
a(a+1)\cdots(a+n-1),\quad \+n\in\mb{N},
\end{cases}
\emn
for any complex number $a$.
Using the fundamental property of gamma function $\Gamma(a+1)=a\Gamma(a)$, $(a)_n$ can be written as
\bmn
(a)_n=\frac{\Gamma(a+n)}{\Gamma(a)}, \quad \quad \quad (n\in\mb{N}\cup \{0\}),
\emn
where $\Gamma$ is the well known Gamma function.

The series ${_p F_q}$ defined by \eqref{def} is converges for all values of $z$ whenever $p\leq q$.
Further, if $p=q+1$, then the series \eqref{def} converges when $|z|<1$. Also it is absolutely convergent on the unite circle $|z|=1$
if $\mathfrak{R}(\sum_{j=1}^q b_j-\sum_{j=1}^p a_j)>0$ and it is convergent on the unit circle $|z|=1$ except at $z=1$
provided $-1<\mathfrak{R}(\sum_{j=1}^q b_j-\sum_{j=1}^p a_j)\leq 0$. For more details about this function, we refer
\cite{kn:rainville, kn:slater}.

On the other hand, we define the Laplace transform of a function $f(t)$ of a real variable $t$ as the integral $g(s)$ over
a range of the complex parameter $s$ by the integral
\bmn\label{def-l}
g(s)=\mathfrak{L}\{f(t); s\}=\int_0^\infty e^{-st}f(t)dt,
\emn
provided this integral exists in the Lebesgue sense. For more details about the Laplace transforms, we refer \cite{kn:dav,kn:doe}.

Now, keeping in mind, the following well known and useful result
\bmn
\int_0^\infty e^{-st}t^{\alpha-1}dt=\Gamma(\alpha)s^{-\alpha},
\emn
provided $\mathfrak{R}(s)>0$ and $\mathfrak{R}(\alpha)>0$.
If we employ \eqref{def} with $p\leq q$, then it is a simple exercise to arrive at the following Laplace transform of
a generalized hypergeometric function ${_p F_q}$ as:
\bmn \label{def-t}
\int_0^\infty e^{-st}t^{v-1}\: {_p F_q}
\ffnk{ccc}{wt}{a_1,\+\cdots,\+a_p}
              {b_1,\+\cdots,\+b_q}dt
=\Gamma(v)s^{-v}{_{p+1} F_q}
\ffnk{cccc}{\frac{w}{s}}{v,\+a_1,\+\cdots,\+a_p}
              {\+b_1,\+\cdots,\+b_q},
\emn
provided $\mathbf{(i)}$ if $p< q$, $\mathfrak{R}(v)>0$, $\mathfrak{R}(s)>0$ and $w$ is arbitrary or
$\mathbf{(ii)}$ if $p=q>0$, $\mathfrak{R}(v)>0$ and $\mathfrak{R}(s)>\mathfrak{R}(w)$,
especially $\mathbf{(iii)}$ if $p=q>0$, $s=w$, $\mathfrak{R}(v)>0$, $\mathfrak{R}(s)>0$ and $\mathfrak{R}(\sum_{j=1}^q b_j-\sum_{j=1}^p a_j-v)>0$.

Further, it is not out of place to mention here that interchanging the order of summation and integration (in the proof of \eqref{def-t})
is easily seen to be justified due to the uniform convergence of the series defined by \eqref{def}.

In particular, when $p=q=1$, for Kummer's confluent hypergeometric function $_1F_1$ (also referred to as the
confluent hypergeometric function of the first kind), we see that its Laplace transform is
\bmn \label{lap-1}
\int_0^\infty e^{-st}t^{b-1}\: {_1 F_1}
\ffnk{c}{wt}{a}{c}dt
=\Gamma(b)s^{-b}{_2F_1}
\ffnk{cc}{\frac{w}{s}}{a,\+b}
              {\+c},
\emn
provided $\mathfrak{R}(b)>0$ and $\mathbf{(i)}$ $\mathfrak{R}(s)>\max\{\mathfrak{R}(w), 0\}$ or
$\mathbf{(ii)}$ $s=w$, $\mathfrak{R}(s)>0$ and $\mathfrak{R}(c-a-b)>0$.

Also, when $p=q=2$ and $p=q=3$ for generalized hypergeometric functions
${_2F_2}$ and ${_3F_3}$, we see their Laplace transforms, respectively, are given by
\bmn \label{lap-2}
\int_0^\infty e^{-st}t^{v-1}\: {_2 F_2}
\ffnk{cc}{wt}{a_1,\+a_2}{b_1,\+b_2}dt
=\Gamma(v)s^{-v}{_3F_2}
\ffnk{ccc}{\frac{w}{s}}{v,\+a_1,\+a_2}
              {\+b_1,\+b_2},
\emn
provided $\mathfrak{R}(v)>0$ and $\mathbf{(i)}$ $\mathfrak{R}(s)>\max\{\mathfrak{R}(w), 0\}$ or
$\mathbf{(ii)}$ $s=w$, $\mathfrak{R}(s)>0$ and $\mathfrak{R}(b_1+b_2-a_1-a_2-v)>0$,

and
\bmn \label{lap-3}
\int_0^\infty e^{-st}t^{v-1}\: {_3 F_3}
\ffnk{ccc}{wt}{a_1,\+a_2,\+a_3}{b_1,\+b_2,\+b_3}dt
=\Gamma(v)s^{-v}{_4F_3}
\ffnk{cccc}{\frac{w}{s}}{v,\+a_1,\+a_2,\+a_3} {\+b_1,\+b_2,\+b_3}.
\emn
provided $\mathfrak{R}(v)>0$ and $\mathbf{(i)}$ $\mathfrak{R}(s)>\max\{\mathfrak{R}(w), 0\}$ or
$\mathbf{(ii)}$ $s=w$, $\mathfrak{R}(s)>0$ and $\mathfrak{R}(b_1+b_2+b_3-a_1-a_2-a_3-v)>0$.

By employing classical summation theorems such as those of Gauss second, Kummer and Bailey for the series $_2F_1$;
Watson, Dixon and Whipple for the series $_3F_2$ and their generalizations
\cite{kn:lavoie-grondin-rathie1,kn:lavoie-grondin-rathie2,kn:lavoie-grondin-rathie3}.
Recently Kim et al. \cite{kn:kim-rathie-cvi1, kn:kim-rathie-cvi2} have obtained a large number of Laplace transforms
for the confluent hypergeometric function $_1F_1$ and generalized hypergeometric function $_2F_2$.
Here, in our present investigation, will mention a few of them, which are:
\bmn\label{gauss2-l}
\int_0^\infty e^{-st}t^{b-1}\:{_1F_1}\ffnk{c}{\frac12 t s}{a}{\frac12(a+b+1)}dt
=s^{-b}\frac{\Gamma(\frac12)\Gamma(b)\Gamma(\frac12a+\frac12b+\frac12)}{\Gamma(\frac12a+\frac12)\Gamma(\frac12b+\frac12)},
\emn
provided $\mathfrak{R}(b)>0$ and $\mathfrak{R}(s)>0$.

\bmn\label{bailey-l}
\int_0^\infty e^{-st}t^{-a}\:{_1F_1}\ffnk{c}{\frac12 t s}{a}{c}dt
=s^{a-1}\frac{\Gamma(1-a)\Gamma(\frac12c)\Gamma(\frac12c+\frac12)}{\Gamma(\frac12a+\frac12c)\Gamma(\frac12c-\frac12a+\frac12)}.
\emn
provided $\mathfrak{R}(1-a)>0$ and $\mathfrak{R}(s)>0$.

\bmn\label{kummer-l}
\int_0^\infty e^{-st}t^{b-1}\:{_1F_1}\ffnk{c}{- t s}{a}{1+a-b}dt
=\frac{s^{-b}2^{-a}\Gamma(\frac12)\Gamma(b)\Gamma(1+a-b)}{\Gamma(\frac12a+\frac12)\Gamma(1+\frac12a-b)},
\emn
provided $\mathfrak{R}(b)>0$ and $\mathfrak{R}(s)>0$.

\bmn\label{watson-l}
\+\+\int_0^\infty e^{-st}t^{c-1}\:{_2F_2}\ffnk{c}{s t}{a, \:\:b}{\frac12(a+b+1), 2c}dt\\[2mm]
\+\+=\frac{s^{-c}\Gamma(\frac12)\Gamma(c)\Gamma(c+\frac12)\Gamma(\frac12a+\frac12b+\frac12)\Gamma(c-\frac12a-\frac12b+\frac12)}
{\Gamma(\frac12a+\frac12)\Gamma(\frac12b+\frac12)\Gamma(c-\frac12a+\frac12)\Gamma(c-\frac12b+\frac12)},\nnm
\emn
provided $\mathfrak{R}(c)>0$, $\mathfrak{R}(s)>0$ and $\mathfrak{R}(2c-a-b)>-1$.

\bmn\label{dixon-l}
\+\+\int_0^\infty e^{-st}t^{c-1}\:{_2F_2}\ffnk{c}{s t}{a, \:\:b}{1+a-b, 1+a-c}dt\\[2mm]
\+\+=\frac{s^{-c}\Gamma(c)\Gamma(1+\frac12a)\Gamma(1+a-b)\Gamma(1+a-c)\Gamma(1+\frac12a-b-c)}
{\Gamma(1+a)\Gamma(1+\frac12a-b)\Gamma(1+\frac12a-c)\Gamma(1+a-b-c)},\nnm
\emn
provided $\mathfrak{R}(c)>0$, $\mathfrak{R}(s)>0$ and $\mathfrak{R}(a-2b-2c)>-2$.

\bmn\label{whipple-l}
\+\+\int_0^\infty e^{-st}t^{c-1}\:{_2F_2}\ffnk{c}{s t}{a, \:\:b}{d, \:\:e}dt\\[2mm]
\+\+=\frac{s^{-c}\pi \Gamma(c)\Gamma(d)\Gamma(e)}
{2^{2c-1}\Gamma(\frac12a+\frac12d)\Gamma(\frac12a+\frac12e)\Gamma(\frac12b+\frac12d)\Gamma(\frac12b+\frac12e)},\nnm
\emn
provided $\mathfrak{R}(c)>0$ and $\mathfrak{R}(s)>0$ with $a+b=1$ and $d+e=1+2c$.

\textbf{Remark.} The results \eqref{gauss2-l} and \eqref{bailey-l} are also recorded in \cito{slater}.

The aim of this research paper is to obtain certain new and useful (potentially) Laplace transforms for the generalized hypergeometric functions
${_2F_2}$ and ${_3F_3}$, so far not recorded in the literature, by using \eqref{lap-2} and \eqref{lap-3}
with the help of known extensions of the classical summation theorems
obtained earlier by Kim et al. \cito{kim-Rakha-rathie1}. For this, we shall require the following summation formulae due to
Kim et al. \cito{kim-Rakha-rathie1}.

\textbf{Extension of Gauss second summation theorem}
\bmn\label{gauss2-e}
\+{_3F_2}\+\ffnk{cc}{\frac12}{a, \quad b, \+d+1}{\frac12(a+b+3),\+d}\\[2mm]
\+\+=\frac{\Gamma(\frac12)\Gamma(\frac12a+\frac12b+\frac32)\Gamma(\frac12a-\frac12b-\frac12)}{\Gamma(\frac12a-\frac12b+\frac32)}
\Big\{\dfrac{\frac12(a+b-1)-\frac{a b}{d}}{\Gamma(\frac12a+\frac12)\Gamma(\frac12b+\frac12)}+
\dfrac{\frac{a+b+1}{d}-2}{\Gamma(\frac12a)\Gamma(\frac12b)}\Big\},\nnm
\emn
provided $\mathfrak{R}(d)>0$.

\textbf{Extension of Bailey summation theorem}

\bmn\label{bailey-e}
\+{_3F_2}\+\ffnk{ccc}{\frac12}{a, \+1-a, \+d+1}{\+c+1, \+d}
=2^{-c}\Gamma(\frac12)\Gamma(c+1)\\[2mm]
\+\+\times\Big\{\frac{\frac{2}{d}}{\Gamma(\frac12a+\frac12c)\Gamma(\frac12c-\frac12a+\frac12)}
+\frac{1-\frac{c}{d}}{\Gamma(\frac12a+\frac12c+\frac12)\Gamma(\frac12c-\frac12a+1)}\Big\},\nnm
\emn
provided $\mathfrak{R}(d)>0$.

\textbf{Extension of Kummer summation theorem}

\bmn\label{kummer-e}
\+{_3F_2}\+\ffnk{cc}{-1}{a, \quad b, \+d+1}{2+a-b,\+d}\\[2mm]
\+\+=\frac{\Gamma(\frac12)\Gamma(2+a-b)}{2^a(1-b)}
\Big\{\frac{\frac{1+a-b}{d}-1}{\Gamma(\frac12a)\Gamma(\frac12a-b+\frac32)}+
\frac{1-\frac{a}{d}}{\Gamma(\frac12a+\frac12)\Gamma(1+\frac12a-b)}\Big\},\nnm
\emn
provided $\mathfrak{R}(d)>0$.

\textbf{First extension of Watson summation theorem}

\bmn\label{watson-e1}
\+{_4F_3}\+\ffnk{cc}{1}{a, \quad\quad b, \quad\quad c \+d+1}{\frac12(a+b+1),\quad 2c+1, \+d}\\[2mm]
\+\+=\frac{2^{a+b-2}\Gamma(c+\frac12)\Gamma(\frac12a+\frac12b+\frac12)\Gamma(c-\frac12a-\frac12b+\frac12)}{\Gamma(\frac12)\Gamma(a)\Gamma(b)}\nnm\\[2mm]
\+\+\times\Big\{\frac{\Gamma(\frac12a)\Gamma(\frac12b)}{\Gamma(c-\frac12a+\frac12)\Gamma(c-\frac12b+\frac12)}
+(\frac{2c-d}{d})\frac{\Gamma(\frac12a+\frac12)\Gamma(\frac12b+\frac12)}{\Gamma(c-\frac12a+1)\Gamma(c-\frac12b+1)}\Big\},\nnm
\emn
provided $\mathfrak{R}(d)>0$ and $\mathfrak{R}(2c-a-b)>-1$.

\textbf{Second extension of Watson summation theorem}

\bmn\label{watson-e2}
\+{_4F_3}\+\ffnk{ccc}{1}{a, \quad b, \+c, \quad \+d+1}{\frac12(a+b+3), \+2c, \+d}\\[2mm]
\+\+=\frac{2^{a+b-2}\Gamma(c+\frac12)\Gamma(\frac12a+\frac12b+\frac32)\Gamma(c-\frac12a-\frac12b-\frac12)}{(a-b-1)(a-b+1)\Gamma(\frac12)\Gamma(a)\Gamma(b)}\nnm\\[2mm]
\+\+\times\Big\{\alpha\cdot\frac{\Gamma(\frac12a)\Gamma(\frac12b)}{\Gamma(c-\frac12a+\frac12)\Gamma(c-\frac12b+\frac12)}
+\beta \cdot\frac{\Gamma(\frac12a+\frac12)\Gamma(\frac12b+\frac12)}{\Gamma(c-\frac12a)\Gamma(c-\frac12b)}\Big\},\nnm
\emn
provided $\mathfrak{R}(d)>0$ and $\mathfrak{R}(2c-a-b)>-1$ with
\[\alpha=a(2c-a)+b(2c-b)-2c+1-\frac{a b}{d}(4c-a-b-1)\quad
\text{and} \quad
\beta=8\{\frac{1}{2d}(a+b+1)-1\}.\]

\textbf{Extension of Dixon summation theorem}

\bmn\label{dixon-e}
\+{_4F_3}\+\ffnk{ccccc}{1}{a, \+ b, \+ c, \+ d+1}{\+ 2+a-b, \+1+a-c, \+d}\\[2mm]
\+\+=\frac{\alpha}{b-1} \cdot\frac{2^{-a}\Gamma(\frac12)\Gamma(2+a-b)\Gamma(1+a-c)\Gamma(\frac12a-b-c+\frac32)}
{\Gamma(\frac12a)\Gamma(2+a-b-c)\Gamma(\frac12a-c+\frac12)\Gamma(\frac12a-b+\frac32)}\nnm \\[2mm]
\+\++\frac{\beta}{b-1} \cdot \frac{2^{-a-1}\Gamma(\frac12)\Gamma(1+a-b)\Gamma(1+a-c)\Gamma(1+\frac12a-b-c)}
{\Gamma(\frac12a+\frac12)\Gamma(1+a-b-c)\Gamma(1+\frac12a-c)\Gamma(1+\frac12a-c)},\nnm
\emn
provided $\mathfrak{R}(d)>0$ and $\mathfrak{R}(a-2b-2c)>-2$ with
\[\alpha=1-\frac{1}{d}(1+a-b)\quad
\text{and} \quad
\beta=\frac{1+a-b}{1+a-b-c}\big\{\frac{a}{d}(1+a-b-2c)-2(1+\frac12a-b-c)\big\}.\]

\textbf{Extension of Whipple summation theorem}

\bmn\label{whipple-e}
\+{_4F_3}\+\ffnk{cccc}{1}{a, \+ 1-a, \+ c, \+d+1}{\+e+1, \+2c-e+1, \+d}\nnm\\[2mm]
\+\+=\frac{2^{-2a}\Gamma(e+1)\Gamma(e-c)\Gamma(2c-e+1)}{\Gamma(e-a+1)\Gamma(e-c+1)\Gamma(2c-a-e+1)}\nnm\\[2mm]
\+\+\times\Big\{(1-\frac{2c-e}{d})\frac{\Gamma(\frac12e-\frac12 a+1)\Gamma(c-\frac12a-\frac12e+\frac12)}{\Gamma(\frac12a+\frac12e)\Gamma(c-\frac12e+\frac12a+\frac12)}\\[2mm]
\+\++(\frac{e}{d}-1) \frac{\Gamma(\frac12e-\frac12 a+\frac12)\Gamma(c-\frac12a-\frac12e+1)}{\Gamma(\frac12a+\frac12e+\frac12)\Gamma(c+\frac12a-\frac12e)}\Big\},\nnm
\emn
provided $\mathfrak{R}(d)>0$ and $\mathfrak{R}(c)>0$.

\section{Three general Laplace transforms of ${_2F_2}(x)$}

In this section, we shall list three general Laplace transforms of ${_2F_2}(x)$ obtained with the help
of \eqref{lap-2} and the extensions of summation formulas \eqref{gauss2-e}, \eqref{bailey-e} and \eqref{kummer-e}.
Clearly, since \eqref{lap-2} is the most general case, so it is desirable to find, as much as possible,
less general case involving various particular values of the parameters $a$, $b$, $d$ and $e$.
Below, in \eqref{gauss2-e-l}, \eqref{bailey-e-l} and \eqref{kummer-e-l}, we give three new and very general Laplace transforms
of ${_2F_2}$, which are not listed in the standard table of the Laplace transforms books \cite{kn:dav,kn:erd,kn:obe-bad,kn:pru}.

\bmn\label{gauss2-e-l}
\+\+\int_0^\infty e^{-st}t^{b-1}\:{_2F_2}\ffnk{cc}{\frac12t s}{a, \+d+1}{\frac12(a+b+3),\+d}dt\\[2mm]
\+\+=\frac{s^{-b}\Gamma(\frac12)\Gamma(b)\Gamma(\frac12a+\frac12b+\frac32)\Gamma(\frac12a-\frac12b-\frac12)}{\Gamma(\frac12a-\frac12b+\frac32)}
\Big\{\frac{\frac12(a+b-1)-\frac{a b}{d}}{\Gamma(\frac12a+\frac12)\:\Gamma(\frac12b+\frac12)}+
\frac{\frac{a+b+1}{d}-2}{\Gamma(\frac12a)\:\Gamma(\frac12b)}\Big\},\nnm
\emn
provided $\mathfrak{R}(b)>0$, $\mathfrak{R}(d)>0$ and $\mathfrak{R}(s)>0$.

\bmn\label{bailey-e-l}
\+\+\int_0^\infty e^{-st}t^{-a}\:{_2F_2}\ffnk{cc}{\frac12t s}{a, \+d+1}{c+1, \+d}dt
=\dfrac{s^{a-1}\Gamma(\frac12)\Gamma(1-a)\Gamma(c+1)}{2^{c}}\nnm\\[2mm]
\+\+\times\Big\{\frac{\frac{2}{d}}{\Gamma(\frac12a+\frac12c)\Gamma(\frac12c-\frac12a+\frac12)}
+\dfrac{1-\frac{c}{d}}{\Gamma(\frac12a+\frac12c+\frac12)\Gamma(\frac12c-\frac12a+1)}\Big\},\nnm
\emn
provided $\mathfrak{R}(1-a)>0$, $\mathfrak{R}(d)>0$ and $\mathfrak{R}(s)>0$.

\bmn\label{kummer-e-l}
\+\+\int_0^\infty e^{-st}t^{b-1}\:{_2F_2}\ffnk{cc}{-t s}{a, \+d+1}{2+a-b,\+d}dt\\[2mm]
\+\+=\frac{s^{-b}\Gamma(\frac12)\Gamma(b)\Gamma(2+a-b)}{2^a(1-b)}
\Big\{\frac{\frac{1+a-b}{d}-1}{\Gamma(\frac12a)\Gamma(\frac12a-b+\frac32)}+
\frac{1-\frac{a}{d}}{\Gamma(\frac12a+\frac12)\Gamma(1+\frac12a-b)}\Big\},\nnm
\emn
provided $\mathfrak{R}(b)>0$, $\mathfrak{R}(d)>0$ and $\mathfrak{R}(s)>0$.

\textbf{Remark:} In \eqref{gauss2-e-l}, \eqref{bailey-e-l} and \eqref{kummer-e-l}, if we respectively take
$d=\frac12(a+b+1)$, $d=c$ and $d=1+a-b$, we recover \eqref{gauss2-l}, \eqref{bailey-l} and \eqref{kummer-l}
obtained earlier by Kim et al. \cito{kim-rathie-cvi1}.

\section{Four general Laplace transforms of ${_3F_3}(x)$}
As described in section 2, here we shall mention from general Laplace transforms of ${_3F_3}(x)$ with the help of \eqref{lap-3}
and the extensions of summation formulas \eqref{watson-e1}, \eqref{watson-e2}, \eqref{dixon-e} and \eqref{whipple-e}.
Below, in \eqref{watson-e1-l}, \eqref{watson-e2-l}, \eqref{dixon-e-l} and \eqref{whipple-e-l}, we give four new and very general Laplace transforms
of ${_3F_3}$, which are not listed in the standard tables of Laplace transforms books \cite{kn:dav,kn:erd,kn:obe-bad,kn:pru}.

\bmn\label{watson-e1-l}
\+\+\int_0^\infty e^{-st}t^{c-1}\:{_3F_3}\ffnk{ccc}{s t}{a, \+ b, \+d+1}{\frac12(a+b+1), \+2c+1, \+d}dt\\
\+\+=\frac{s^{-c}2^{a+b-2}\Gamma(c)\Gamma(c+\frac12)\Gamma(\frac12a+\frac12b+\frac12)\Gamma(c-\frac12a-\frac12b+\frac12)}{\Gamma(\frac12)\Gamma(a)\Gamma(b)}\nnm\\
\+\+\times\Big\{\frac{\Gamma(\frac12a)\Gamma(\frac12b)}{\Gamma(c-\frac12a+\frac12)\Gamma(c-\frac12b+\frac12)}
+(\frac{2c-d}{d})\frac{\Gamma(\frac12a+\frac12)\Gamma(\frac12b+\frac12)}{\Gamma(c-\frac12a+1)\Gamma(c-\frac12b+1)}\Big\},\nnm
\emn
provided $\mathfrak{R}(c)>0$, $\mathfrak{R}(s)>0$ and $\mathfrak{R}(2c-a-b)>-1$.

\bmn\label{watson-e2-l}
\+\+\int_0^\infty e^{-st}t^{c-1}\:{_3F_3}\ffnk{ccc}{s t}{a, \+ b,\+d+1}{\frac12(a+b+3), \+2c, \+d}dt\\[2mm]
\+\+=\frac{s^{-c}2^{a+b-2}\Gamma(c)\Gamma(c+\frac12)\Gamma(\frac12a+\frac12b+\frac32)\Gamma(c-\frac12a-\frac12b-\frac12)}
{(a-b-1)(a-b+1)\Gamma(\frac12)\Gamma(a)\Gamma(b)}\nnm\\[2mm]
\+\+\times\Big\{\alpha\cdot\frac{\Gamma(\frac12a)\Gamma(\frac12b)}{\Gamma(c-\frac12a+\frac12)\Gamma(c-\frac12b+\frac12)}
+\beta \cdot\frac{\Gamma(\frac12a+\frac12)\Gamma(\frac12b+\frac12)}{\Gamma(c-\frac12a)\Gamma(c-\frac12b)}\Big\},\nnm
\emn
provided $\mathfrak{R}(c)>0$, $\mathfrak{R}(s)>0$ and $\mathfrak{R}(2c-a-b)>-1$ with
\[\alpha=a(2c-a)+b(2c-b)-2c+1-\frac{a b}{d}(4c-a-b-1)\quad
\text{and} \quad
\beta=8\{\frac{1}{2d}(a+b+1)-1\}.\]

\bmn\label{dixon-e-l}
\+\+\int_0^\infty e^{-st}t^{c-1}\:{_3F_3}\ffnk{ccc}{s t}{a, \+ b, \+ d+1}{2+a-b, \+1+a-c,\+ d}dt\\[2mm]
\+\+=\frac{s^{-c}2^{-a}\Gamma(\frac12)\Gamma(c)}{b-1}
\Big\{\alpha\cdot\frac{\Gamma(2+a-b)\Gamma(1+a-c)\Gamma(\frac12a-b-c+\frac32)}
{\Gamma(\frac12a)\Gamma(2+a-b-c)\Gamma(\frac12a-c+\frac12)\Gamma(\frac12a-b+\frac32)}\nnm \\[2mm]
\+\++\frac{\beta}{2} \cdot \frac{\Gamma(1+a-b)\Gamma(1+a-c)\Gamma(1+\frac12a-b-c)}
{\Gamma(\frac12a+\frac12)\Gamma(1+a-b-c)\Gamma(1+\frac12a-b)\Gamma(1+\frac12a-c)}\Big\},\nnm
\emn
provided $\mathfrak{R}(c)>0$, $\mathfrak{R}(s)>0$ and $\mathfrak{R}(a-2b-2c)>-2$ with
\[\alpha=1-\frac{1}{d}(1+a-b)\quad
\text{and} \quad
\beta=\frac{1+a-b}{1+a-b-c}\big\{\frac{a}{d}(1+a-b-2c)-2(1+\frac12a-b-c)\big\}.\]

\bmn\label{whipple-e-l}
\+\+\int_0^\infty e^{-st}t^{c-1}\:{_3F_3}\ffnk{ccc}{s t}{a, \+ 1-a, \+d+1}{e+1, \+2c-e+1, \+d}\\[2mm]
\+\+=\frac{s^{-c}2^{-2a}\Gamma(c)\Gamma(e+1)\Gamma(e-c)\Gamma(2c-e+1)}{\Gamma(e-a+1)\Gamma(e-c+1)\Gamma(2c-a-e+1)}\nnm\\[2mm]
\+\+\times\Big\{(1-\frac{2c-e}{d})\frac{\Gamma(\frac12e-\frac12 a+1)\Gamma(c-\frac12a-\frac12e+\frac12)}{\Gamma(\frac12a+\frac12e)\Gamma(c-\frac12e+\frac12a+\frac12)}\nnm\\[2mm]
\+\++(\frac{e}{d}-1) \frac{\Gamma(\frac12e-\frac12 a+\frac12)\Gamma(c-\frac12a-\frac12e+1)}{\Gamma(\frac12a+\frac12e+\frac12)\Gamma(c-\frac12e+\frac12a)}\Big\},\nnm
\emn
provided $\mathfrak{R}(c)>0$ and $\mathfrak{R}(s)>0$.

\textbf{Special cases:}
In \eqref{watson-e1-l}, if we take $d=2c$ or in \eqref{watson-e2-l}, if we take $d=\frac12(a+b+1)$, we recover \eqref{watson-l}.
While in \eqref{dixon-e-l} and \eqref{whipple-e-l}, if we take $d=1+a-b$ and $d=e$, we respectively recovered
\eqref{dixon-l} and \eqref{whipple-l}.

\textbf{Concluding remark:}
we conclude this paper by remarking that the Laplace transforms for the generalized hypergeometric functions
${_2F_2}$ and ${_3F_3}$ established in this paper may be useful in Mathematics, Statistics, Physics and Engineering.

\section*{Acknowledgement}
This work was, in part, supported by National Natural Science Foundation of China (Grant No. 11201291),
Natural Science Foundation of Shanghai (Grant No. 12ZR1443800) and
a grant of "The First-class Discipline of Universities in Shanghai".


\end{document}